\documentclass[10pt]{article}
\font\smallit = cmti10

\usepackage{amssymb,latexsym,amsmath,epsfig,amsthm} 

\usepackage{amsfonts}
\usepackage{graphicx}
\usepackage{float}
\usepackage{hyperref}
\usepackage{placeins}

\makeatletter

\renewcommand\section{\@startsection {section}{1}{\z@}
{-30pt \@plus -1ex \@minus -.2ex}
{2.3ex \@plus.2ex}
{\normalfont\normalsize\bfseries\boldmath}}

\renewcommand\subsection{\@startsection{subsection}{2}{\z@}
{-3.25ex\@plus -1ex \@minus -.2ex}
{1.5ex \@plus .2ex}
{\normalfont\normalsize\bfseries\boldmath}}

\renewcommand{\@seccntformat}[1]{\csname the#1\endcsname. }

\makeatother

\newtheorem{theorem}{Theorem}

\theoremstyle{definition}

\begin{document}


\begin{center}
\uppercase{\bf \boldmath A direct approach to\\the Gallai-Schur numbers}
\vskip 20pt
{\bf Fred Rowley}\\
{\smallit formerly of Lincoln College, Oxford, UK.}\\
{\tt fred.rowley@ozemail.com.au}\\ 
\end{center}
\vskip 10pt

\vskip 30pt


\centerline{\bf Abstract}

\noindent


This paper characterises the structure of every maximal weak or strong Gallai-Schur partition.  

The results confirm the exact values of Gallai-Schur numbers provided by Budden (2020) in the strong case, and provide corresponding values for weak Gallai-Schur numbers.  The proofs are elementary and standalone.  

\pagestyle{myheadings}
\markright{A Direct Approach to the Gallai-Schur Numbers (ArXiv version 1)\hfill}
\thispagestyle{empty}
\baselineskip=12.875pt
\vskip 30pt


\section{Introduction}

In \cite{Bud}, Budden deduced the exact values of the Gallai-Schur numbers (in the strong case, as defined here) from an earlier result of Chung and Graham \cite{ChuGra}.  

These numbers represent upper bounds on the size of colourings (or partitions) of integer intervals $[1, n]$ which avoid both monochromatic and rainbow sums of the form $a+b=c$. We distinguish the \textit{strong case}, where we seek to avoid monochromatic sums where $a=b$, from the \textit{weak case}, where such sums are permitted.  

The result of Chung and Graham in \cite{ChuGra} is graph-theoretic, and supported Budden's approach in the strong case very well, but it does not appear to lend itself easily to the weak case.  Budden mentioned the the weak case, but left its investigation for later work.  Robertson \textit{et al} (2024) derived upper and lower bounds for the weak numbers \cite{Robetal}

In this paper, we provide a self-contained proof of Budden's result, and a parallel result for the weak case.  We prove that the weak Gallai-Schur numbers are in fact equal to the lower bound derived in \cite{Robetal}.

\section{Terminology and Definitions}

In concept, the Gallai-Schur partition follows on from the Schur partition, so we begin there.

If an integer interval $U = [1, n]$ can be partitioned into $r$ disjoint subsets $S_i$ for $i = 1, 2, \dots, r$, where no subset contains any set of integers $a, b, c$, such that $a+b=c$, then each such subset is said to be {\it sum-free} and that partition is called a {\it Schur (r-)partition}. The order of the set $U$ is clearly $n$, and is also referred to as the order of the partition. We may occasionally refer to this type of partition as a {\it strong} partition, to distinguish it from the weak case defined next below.

If we limit the definition above so that it excludes only sets of three {\bf distinct} integers $a, b, c$, such that $a+b=c$, then each such subset is \textit{weakly} sum-free and that partition is a \textit{weak Schur (r-)partition}.  We define a \textit{weak pair} of integers as a pair of the form ${a, 2a}$.  The strong partition excludes such pairs from each subset, whereas the weak partition does not.

A Schur partition is a \textit{Gallai-Schur partition}, if in addition there is no triple $a, b, c$, such that $a+b=c$, and each of $a, b, c$ is in a different subset, noting that the partition must contain exactly $r$ non-empty subsets.  We define a \textit{weak Gallai-Schur partition} correspondingly.  

For any $r$, the \textit{Gallai-Schur number} $GS(r)$ is the least integer $n$ such that no Gallai-Schur r-partition of $[1, m]$ exists for any\footnote{This definition may seem slightly longer than necessary, but it must deal with the fact that choosing $m = r$ forces at least ${1, 2, 3}$ to be a rainbow triple.} integer $m \ge n$.  We may define the \textit{weak Gallai-Schur number} $WGS(r)$ in an exactly parallel manner.  

In section 3, we define some mappings that we will use in section 4 to determine the structure of maximal (strong or weak) Gallai-Schur partitions, and in section 5 we deduce $GS(r)$ and $WGS(r)$ directly. 

In section 6, some brief conclusions are drawn.

\section{Gallai-Schur Partitions}

We begin by defining two very simple constructions, which we will use iteratively.  These are derived by generalising and extending the construction used by Budden \cite{Bud} to demonstrate lower bounds. They are used here, along with their inverses, to derive both upper and lower bounds.   

\begin{theorem}
  \label{Thm:Construx-Thm1}
(2-fold / 5-fold Construction)\\
If there is a (weak) Gallai-Schur partition of the integers $[1, m]$ into $r$ non-empty subsets, then there is a (weak) Gallai-Schur  partition of $[1, 2m+1]$ into $r+1$ non-empty subsets; and a (weak) Gallai-Schur partition of $[1, 5m+4]$ into $r+2$ non-empty subsets.    
\end{theorem}

\begin{proof}

Let $P_r$ be such a partition of $[1, m]$.  We can identify the colour of an integer with the index number of the unique subset $S_i \subset P_r$ that it belongs to. Thus for an integer $a \in S_i$, we write the colour as $\chi(a) = i$.

The first part of the theorem can be demonstrated by considering the partition $Q_{r+1}$ of $[1, 2m+1]$ that results if we assign all the odd numbers to a single new subset $S'_1$.  The members of each of the other new subsets $S'_{t+1} \subset Q_{r+1}$ can be determined by multiplying the members of $S_t \subset P_r$ by 2, for $1 \le t \le r$.  There is no 'weak pair' $(a, 2a)$ in $S'_1$, so we can deduce that if $P_r$ is a strong or weak partition, $Q_{r+1}$ has the same property.

In this way, for example, a strong partition $P_2$ of $[1,  4]$ into $S_{1} = \{1, 4\}$, $S_{2} = \{2, 3\}$ gives rise to $Q_3$ with $S'_1 = \{1, 3, 5, 7, 9\}, S'_2 = \{2, 8\}, S'_3 = \{4, 6\}$ and $2m+1 = 9$.  The subsets are clearly distinct and fully cover the interval $[1, 9]$.  There is no 'weak pair' $(a, 2a)$ in $S'_1$, so $Q_3$ is strong, like $P_2$.

For any triple with $a + b = c$, either (i) exactly two of them are odd, or (ii) none of them is odd.  In case (i), it is easy to check that the sum or difference of two members of $S'_1$ cannot be of the same colour, since the sum or difference is even: and of course $a, b, c$ is not a rainbow triple if two of them are members of $S'_1$.  In case (ii), if there is either a monochromatic or a rainbow triple in the other subsets, then we would know there was such a triple in the partition $P_r$, which is a contradiction.

The mapping from $P_r$ to $Q_{r+1}$ is clearly well-defined, and has an obvious inverse.  We may call the mapping $^2\Theta$.  This proves the first part of the theorem. 

\medskip

The second part can be demonstrated by assigning all integers $x$ in the range $[1, 5m+4]$ to a subset $S''_1$ of a new partition $R_{r+2}$ as follows.
If $x \equiv  2$ or 3 \textit{modulo} 5, then $x \in S''_1$.  If $x \equiv$ 1 or 4 \textit{modulo} 5, then $x \in S''_2$.
The members of each of the other subsets $S''_{t+2}$ can be determined by multiplying the members of $S_t$ by 5, for $1 \le t \le r$.  Once more, the subsets are clearly distinct and fully cover the interval $[1, 5m+4]$.  There is no 'weak pair' $(a, 2a)$ in $S''_1$ or $S''_2$, so if $P_r$ is a strong or weak partition, $R_{r+2}$ has the same property.

It is easy to see that the sum of two members of $S''_1$ must be in $S''_2$, and \textit{vice versa}, and the resulting sum is obviously not rainbow in either case.  Adding a member of $S''_1$ to a member of $S''_2$ clearly cannot form a monochromatic triple, but that triple is never rainbow either, since the sum cannot be a multiple of 5.  Thus we have dealt with the case where no member of the triple is divisible by 5.

If any member of a triple with $a + b = c$ is divisible by 5, then either (i) all of them are divisible by 5, or (ii) exactly one of them is. 

In case (i), we note that if there is either a monochromatic or a rainbow triple contained wholly in the subsets with index greater than 2, then we would know there was such a triple in the partition $P_r$, which is a contradiction.

In case (ii), we must have two members of either $S''_1$ or $S''_2$ whose absolute difference is a multiple of 5, and hence they are in the same subset.  Thus, the triple is neither monochromatic nor rainbow.  

The mapping from $P_r$ to $R_{r+2}$ is clearly well-defined, and has an obvious inverse.  We may call the mapping $^5\Theta$.  This proves the second part of the theorem.  

\end{proof}

\section{Structure of Maximal Partitions}
	
The result in the previous section provides a basis for lower bounds on $GS(r)$ and $WGS(r)$.  In this section, we explore the structure of {\bf maximal} examples of Gallai-Schur partitions. 

The assumption that a partition is maximal -- i.e. that there is no larger Gallai-Schur partition into the same number of colours -- is remarkably powerful in this instance.  

We aim to prove the theorem below.  To do so, we assume here that the colours of any partition occur in the natural order -- so that if we have two colours $i, j$, and $i < j$, then the smallest member of subset $S_i$ is less than the smallest member of $S_j$.  We note that the main arguments in this section apply with the same force whether the strong or the weak criterion applies.  

\begin{theorem}
  \label{Thm:Inverse-Thm2}
(2-fold / 5-fold Image)\\
Let the (strong or weak) Gallai-Schur partition $P_r$ be a maximal partition of order $m$, with $r > 3$.  The first three integers must have colours 1, 2, 2 or 1, 2, 1 respectively.  In the first case, $m$ has residue 4 \textit{modulo} 5 and $P_r$ is an image, under the mapping $^5\Theta$, of a Gallai-Schur partition (strong or weak, respectively) of order $(m-4)/5$.  In the second case, $m$ is odd, and $P_r$ is an image, under the mapping $^2\Theta$, of a Gallai-Schur partition (strong or weak, respectively) of order $(m-1)/2$.   
\end{theorem}

\begin{proof}

Let $R_r$ be a partition of $[1, m]$ into $r$ colours, where $r > 3$.  We assume there is no larger Gallai-Schur partition into $r$ colours, and therefore say that it is \textit{maximal}.  

If we look at the colours of $R_r$ in increasing sequence, we can derive some very useful results.

We recall first that the colour of an integer $i$ is written $\chi(i)$.  So by our assumption above, $\chi(1) = 1$.

We must consider whether we can have $\chi(2) = 1$.  If that were so, then with $r > 3$, there would be a contradiction, since we know that for some $s > 2$ we would have $\chi(s) = 2$, and for some $t > s$ we must have $\chi(t) = 3$. We can assume without loss of generality that these are the smallest integers of their respective colours. We can immediately deduce that $\chi(t-1) = \chi(t-2) = 1$, since otherwise we would have rainbow triples $(1, t-1, t)$ and $(2, t-2, t)$.  However, we now have a monochromatic triple $(1, t-1, t)$.  Thus we must take $\chi(2) = 2$.

It is worth pausing to note that the reasoning above ensures a full reflection symmetry.  That is, if $t > 2$ is the smallest integer of a given colour, then all smaller integers must conform to a symmetrical colour pattern within the interval $[1, t-1]$, in order to avoid a rainbow sum.  

If $\chi(1) = 1$ and $\chi(2) = 2$, then we must have $\chi(3) \in \{1, 2\}$, to avoid a rainbow colouring.  

We note here that since $r > 3$, we cannot have $\chi(4) = 2$ by the same symmetry argument as above.  That is, for some minimal $t > 4$ we would have $\chi(t) = 3$ and therefore $\chi(t-2) = \chi(t-4) = 2$, which avoids a rainbow sum, but must produce a monochromatic triple $(2, t-4, t-2)$.  We note that this argument applies in both weak and strong cases.

We go back to $\chi(3)$.  If we first assume $\chi(3) = 2$, then $\chi(4) \in \{1, 2\}$ to avoid a rainbow colouring, but we know (from the above) that this means $\chi(4) = 1$.  We now have the   pattern ``1221" for the first four integers, and we can see this implies $\chi(5) \notin \{1, 2\}$ if we are to avoid a monochromatic sum. Setting $\chi(5) = 3$ does not produce a rainbow colouring, since the preceding sequence is symmetrical. So we now have ``12213". 

The fact that $\chi(5) = 3$ has interesting consequences for values of $s \in [6, 9]$.  To avoid a rainbow sum, we can deduce that $\chi(s) \in \{3, \chi(s-5)\}$. However, we also know that $\chi(6) \in \{\chi(2), \chi(4)\}$, and $\chi(7) \in \{\chi(3), \chi(4)\}$. Thus $\chi(6) = 1$, and $\chi(7) = 2$. Similarly, we deduce that $\chi(8) = 2$, and $\chi(9) = 1$. It is easy to see that $\chi(10) \notin {1, 2}$.  This pattern of colouring continues indefinitely for the integers that are not multiples of 5.  The fact that no multiple of 5 can be of colour 1 or 2, forces repetition of the ``1221[$\chi(5k)$]" pattern iteratively for every subsequent set of 5 consecutive integers. 

This can easily be formalised inductively for all integers in $R_r$. We assume that, for $t \ge 0$, every integer of the form $5s$ with $1 \le s \le t$ is not of colour 1 or 2; that any integer of the form $5t+1$ or $5t+4$ is colour 1; and that any integer of the form $5t+2$ or $5t+3$ is of colour 2. We have proven this for $t \le 1$, so now consider $t > 1$.\\

Since $\chi(6) = \chi(5t-1) = 1$ we deduce that $\chi(5(t+1)) \ne 1$.

Since $\chi(7) = \chi(5t-2) = 2$ we deduce that $\chi(5(t+1)) \ne 2$. \\

Let $\chi(5t) = x$.  The absence of rainbow triples give us the following:

$\chi(6) = 1$ implies $\chi(5t+6) \in \{1, x\}$. $\chi(7) = 2$ implies $\chi(5t+7) \in \{2, x\}$.

$\chi(8) = 2$ implies $\chi(5t+8) \in \{2, x\}$. $\chi(9) = 1$ implies $\chi(5t+9) \in \{1, x\}$.

\vskip 6pt

However, we can also deduce that:

$\chi(4) = 1$ and $\chi(5t+2) = 2$ implies $\chi(5t+6) \in \{1, 2\}$, so $\chi(5t+6) = 1$.

$\chi(4) = 1$ and $\chi(5t+3) = 2$ implies $\chi(5t+7) \in \{1, 2\}$, so $\chi(5t+7) = 2$.

$\chi(7) = 2$ and $\chi(5t+1) = 1$ implies $\chi(5t+8) \in \{1, 2\}$, so $\chi(5t+8) = 2$.

$\chi(8) = 2$ and $\chi(5t+1) = 1$ implies $\chi(5t+9) \in \{1, 2\}$, so $\chi(5t+9) = 1$.

\newpage

This completes the induction.  Summarising, for all $t \ge 0$,

$\chi(5t+5) = \chi(5(t+1)) \notin \{1, 2\}$.

$\chi(5(t+1)+1) = \chi(5(t+1)+4) = 1$, and  $\chi(5(t+1)+2) = \chi(5(t+1)+3) = 2$.

\vskip 6pt

We see that the integers in $R_r$ that are divisible by 5 are not of colour 1 or 2.  Those integers must therefore be coloured in such a way that if they are all divided by 5, the resulting sets will themselves form a strong or weak Gallai-Schur partition, according to the same property of $R_r$, using two fewer colours than $R_r$.  

We now also see that $m$ must have residue 4 \textit{modulo 5}: if not, we could extend the interval by adding members to $S_1$ or $S_2$ to make it so -- which would contradict the assumed maximality of $R_r$.  (It should be clear from the above that the extension of this pattern of colouring integers not divisible by 5 can never give rise to monochromatic or rainbow sums.)

Summarising the above, the members of $S_1$ consist exactly of all the members of $[1, m]$ with residue 1 or 4 \textit{modulo} 5, and the members of $S_2$ consist exactly of all the members of $[1, m]$ with residue 2 or 3 \textit{modulo} 5.  All members of all the other subsets must be divisible by 5.

The multiples of 5 in that pattern form a subset of the interval $[1, m]$, such that if any two members of that subset form part of an additive triple within the interval, then the other member of that triple is a member of the same subset.  

Thus we can form a new partition of $[1, (m-4)/5]$ by taking the subsets $S_t$ of $R_r$ for $t > 2$, dividing their members by 5, and re-labelling them with indices reduced by 2.  We have clearly created the inverse of the mapping $^5\Theta$.  Since there are no weak pairs in $S_1$ or $S_2$, the new partition is strong or weak according to the status of $R_r$.  This proves the first part of the theorem.

Moving to the second part, the only alternative assumption for $\chi(3)$ to avoid a rainbow colouring is that $\chi(3) = 1$.  Since $r > 3$, by the argument of reflection used above, we cannot have $\chi(4) = 2$.  To avoid a monochromatic sum, we cannot have $\chi(4) = 1$.  Therefore $\chi(4) = 3$.  

For $\chi(5)$, we must look at the two possible triples with $c = 5$.  We deduce that $\chi(5) \in \{1, 3\}$ and  $\chi(5) \in \{1, 2\}$ to avoid rainbow sums.  Thus  $\chi(5) = 1$, and we have the initial colour sequence ``12131''.  We can now see that no even integer up to 8 can be of colour 1, since we avoid monochromatic sums.  We also know $\chi(6)$ must be either 2 or 3, to avoid a rainbow sum.  Whichever it is, we know that $\chi(7) \in \{1, 3\}$ and  $\chi(7) \in \{1, 2\}$.  Thus  $\chi(7) = 1$.  This pattern can also be seen to repeat indefinitely: all the odd integers are members of $S_1$ and all the even integers are in other subsets.  A formal proof is omitted, although only a simpler version of the proof above is needed.  We observe that the order $m$ must be odd, since otherwise we could add the odd integer $(m+1)$ to $S_1$, which would be a contradiction of maximality.

\newpage

As in the first part, we can form a new partition of $[1, (m-1)/2]$ by taking the subsets $S_t$ of $R_r$ for $t > 1$, dividing their members by 2, and re-labelling them with indices reduced by 1.  We have clearly created the inverse of the mapping $^2\Theta$.  Since there are no weak pairs in $S_1$, the new partition is again strong or weak according to the status of $R_r$.  This proves the second part of the theorem, which now includes all possible colourings of a maximal partition.

\end{proof}

\section{Size of Maximal Partitions}
	
To examine the size of a maximal Gallai-Schur partition, we start with any such partition $Q_r$ of $[1, m]$ into exactly $r$ colrs. To better manage the arithmetic, we define a simple function $g(Q_r) = | Q_r | + 1= (m+1)$.  

Theorem \ref{Thm:Inverse-Thm2} tells us that provided $r > 3$, $Q_r$ is the image of another Gallai-Schur partition under either $^5\Theta$ or $^2\Theta$, and that we can determine which is the mapping concerned and apply the inverse mapping. If we apply the inverse of $^5\Theta$, we will get a new partition $Q_{r-2}$ of order $(m-4)/5$ and we therefore have $g(Q_{r-2}) = g(Q_r)/5$. If we apply the inverse of $^2\Theta$, we will get a new partition $Q_{r-1}$ of order $(m-1)/2$ and we therefore have $g(Q_{r-1}) = g(Q_r)/2$. 

We can continue this process until the number of colours remaining falls below 4. We will then have a partition $Q_{r-t}$ where $r - t = 2$ or 3.

At that point we may consider how many times we have applied $^2\Theta^{-1}$. We find it cannot be more than once.  If it were more, we could take $Q_{r-t}$ and apply a different sequence of $^i\Theta$ (replacing two instances of $^2\Theta$ by one instance of $^5\Theta$) to generate a $Q'_r$ with $g(Q'_{r}) > g(Q_r)$, which is a contradiction on maximality.

Clearly, the partition $Q_{r-t}$ must also be maximal, otherwise we could easily generate a larger $Q'_r$.

The maximal Gallai-Schur partitions for $s$ colours, $s < 4$ are highly constrained and are easily proven to be as follows.  We use two forms of notation which are equivalent, the first specifying subsets and the second listing colours of the integers in increasing order.

\medskip

\textbf{Strong case:}

\medskip
$B_1 = \{1\} \equiv$ ``1''. $|B_1| = 1.$\newline

$B_2 = \{1, 4\},~\{2, 3\} \equiv$ ``1221''. $|B_2| = 4.$\newline

$B_{3A} = B_1*{^5\Theta} = \{1, 4, 6, 9\}, \{2, 3, 7, 8\}, \{5\} \equiv$ ``122131221''. $|B_{3A}| = 9.$ \newline

$B_{3B} = B_2*{^2\Theta} = \{1, 3, 5, 7, 9\}, \{2, 8\}, \{4, 6\} \equiv$ ``121313121''. $|B_{3B}| = 9.$ \newline

\medskip
\medskip
\medskip

\textbf{Weak case:}

\medskip
$C_1 = \{1, 2\} \equiv$ ``11''. $|C_1| = 2.$\\

$C_2 = \{1, 2, 4, 8\}, \{3, 5, 6, 7\} \equiv$ ``11212221''. $|C_2| = 8.$\\

$C_3 = C_2*{^2\Theta} = \{1, 3, 5, 7, 9, 11, 13, 15, 17\}, \{2, 4, 8, 16\}, \{6, 10, 12, 14\}$ \\ $\equiv$ ``12121312131313121''. $|C_3| = 17.$\\
\\
To be certain, these have been verified as being the only maximal cases, down to re-labelling, using a SAT-solver (Penelope).\\
\\
Returning now to the main argument, we have several sub-cases to consider.

\medskip
\textbf{Strong case:}

If $r - t = 2$ we must have $Q_{r-t} = B_2$.  If $r$ is even we therefore have $g(Q_r) = g(B_2).5^{(r-2)/2} = 5.(5^{(r-2)/2})$.  If $r$ is odd, we must have $g(Q_r) = 2.g(B_2).5^{(r-3)/2} = 10.(5^{(r-3)/2})$.

If $r - t = 3$ we must have $Q_{r-t} = B_{3A}$ or $B_{3B}$.  If $r$ is odd, we therefore have $g(Q_{r-t}) = g(B_{3A}) = g(B_{3B}) = 10$ and $g(Q_r) = 10.(5^{(r-3)/2})$ as above.  If $r$ is even, we can produce a larger partition $Q'_r$ by setting $Q_{r-t} = B_{3B}$ and replacing two applications of $^2\Theta$ with one application of $^5\Theta$, showing that $Q_r$ is not maximal.  This contradiction shows we cannot have $r - t = 3$ in a maximal partition when $r$ is even and $Q_r$ is maximal.

Thus in the strong case, for $r \ge 1$ we have:\\
\indent \indent \indent $g(Q_r) = 5^{(r/2)}$ for even $r$, and \\
\indent \indent \indent $g(Q_r) = 2.(5^{(r-1)/2})$ for odd $r$. 

\medskip
\textbf{Weak case:}

If $r - t = 2$, we must have $Q_{r-t} = C_2$. Since we know there cannot be more than one occurrence of $^2\Theta^{-1}$ in the sequence of inverse mappings applied, we must have $g(Q_r) = g(C_2).5^{(r-2)/2} = 9.(5^{(r-2)/2)})$ if $r$ is even. If $r$ is odd, $g(Q_r) = 2.g(C_2).5^{(r-3)/2} = 18.(5^{(r-3)/2)})$.

If $r - t = 3$, we must have $Q_{r-t} = C_3$. Since we know there cannot be more than one occurrence of $^2\Theta^{-1}$ in the sequence of inverse mappings applied, we must again have $g(Q_r) = 2.g(C_3).5^{(r-2)/2} = 18.(5^{(r-3)/2})$ if $r$ is odd.  If $r$ is even, we have one or more occurrences of $^2\Theta^{-1}$ in the sequence of inverse mappings, so the value of $g(Q_r)$ will be less than $g(C_2).5^{(r-2)/2}$, hence $Q_r$ is not maximal. This contradiction shows we cannot have $r - t = 3$ in the weak case when $r$ is even and $Q_r$ is maximal.

\newpage

It should be clear that since each of the $Q_r$ is maximal, $g(Q_r)$
corresponds to the (strong or weak) Gallai-Schur number.  Thus we have proved the following:

\begin{theorem}
  \label{Thm:S&W G-S Numbers}
(Strong and weak Gallai-Schur numbers)\\
In the strong case, for $r \ge 1$ we have:\\
\indent \indent \indent $GS(r) = 5^{(r/2)}$ for even $r$, and \\
\indent \indent \indent $GS(r) = 2.(5^{(r-1)/2})$ for odd $r$.\\
In the weak case, we have the exceptional value
$WGS(r) = 3$ for $r = 1$,\\
and for $r > 1$,\\
\indent \indent \indent $WGS(r) = 9.(5^{(r-2)/2})$ for even $r$, and \\
\indent \indent \indent $WGS(r) = 18.(5^{(r-3)/2})$ for odd $r$.

\vskip 6pt
\indent \indent \indent i.e. $WGS(r) = (9/5).GS(r)$ for $r > 1$. 
\end{theorem}

\section{Conclusions}
	
We have precisely characterised maximal Gallai-Schur partitions, and determined the Gallai-Schur numbers.  It is noted that both the weak and strong maximal Gallai-Schur partitions are unique for even $r$, although not for odd $r$.

This is possible because the combination of constraints that require a triple to be free of both rainbow and monochromatic sums is very strong, so that maximal partitions take a restricted range of forms.  This in turn allows us to derive exact values for weak and strong Gallai-Schur numbers, replicating, for the strong case, those derived by Budden in \cite{Bud}. 

\vskip 12pt

\noindent {\bf Dedication}

I dedicate this paper to the memory of my late father, Gordon Rowley, formerly of Queens' College, Cambridge, who took on some memorable and fulfilling challenges, late in life.


\medskip
















\end{document}